\documentclass[12pt]{article}
\usepackage{latexsym}
\usepackage[dvips]{graphicx}
\newcommand{\dstyle}[1]{\displaystyle{#1}}
\newcommand{\defeq}{\stackrel{\rm def}{=}}
\newcommand{\e}{{\rm e}}

\newcommand{\del}[1]{\Delta #1}
\newcommand{\bsquare}{\hbox{\rule{6pt}{6pt}}}
\renewcommand{\d}{{\rm d}}
\newtheorem{Lem_eng}{Lemma}[section]
\makeatletter

 \@addtoreset{equation}{section}
\makeatother
\begin{document}
\begin{center}
  {\LARGE
    Properties of the Weibull cumulative exposure model
  }
\end{center}
\begin{center}
YOSHIO KOMORI\\
Department of Systems Innovation and Informatics,\\
Kyushu Institute of Technology,
Iizuka 820-8502, Japan
\end{center}
%

This article is aimed at the investigation of some properties of
the Weibull cumulative exposure model on multiple-step step-stress
accelerated life test data.
Although the model includes a probabilistic idea of Miner's rule
in order to express the effect of cumulative damage in fatigue,
our result shows that the application of only this is
not sufficient to express degradation
of specimens and the shape parameter must be larger than 1.
For a random variable obeying the model, its average and standard
deviation are investigated on a various sets of parameter values.
In addition,
a way of checking the validity of the model is illustrated
through an example of the maximum likelihood estimation
on an actual data set, which is about time to breakdown
of cross-linked polyethylene-insulated cables.
%
\vskip 4mm
%
\section{Introduction}

In many industrial fields it is requested for lots of products to
operate for a long period of time. Accompanied with that, it is very
important to give reliability in relation to the lifetime of products.
In such cases, however, life testing under a normal stress can lead
to a lengthy procedure with expensive cost. As a means to cope with
these problems, the study of accelerated life test (ALT) has been
developed. The test makes it possible to quickly obtain information
on the life distribution of products by inducing early failure with
stronger stress than normal.

One important way in ALT is step-stress accelerated life test (SSALT).
There are mainly two types of SSALTs, a simple SSALT and a multiple-step
SSALT. In the simple SSALT there is a single change of stress during the
test.
Miller and Nelson
(1983)
have shown optimum simple SSALT plans
in an exponential cumulative exposure (CE) model.
Xiong
(1999)
has studied an exponential CE model
with a threshold parameter in the simple SSALT.
Park and Yum
(1998)
have shown optimum modified
simple SSALT plans in an exponential CE model,
under the consideration that it is desirable to
increase the stress at some finite rate.
Lu and Rudy
(2002)
have dealt with the Weibull CE
model with the inverse power law
in the simple SSALT.

On the other hand, in the multiple-step SSALT there are changes of stress
more than once. Yeo and Tang
(1999)
have investigated a
three-step SSALT
in an exponential CE model. Khamis
(1997)
has proposed
an exponential CE model with $k$ explanatory variables and
investigated it on three-step SSALT data. McSorley, Lu and
Li
(2002)
have shown the properties
of the maximum likelihood (ML) estimators
of parameters in the Weibull CE model
with a log-linear function of stress
on three-step SSALT data.
Nelson
(1980, 1990)
has proposed an important idea, which
gives the basic CE model for life as a function of
constant stress from SSALT data.
This is a probabilistic analog of Miner's rule
(Miner, 1945),
which is stated on a deterministic situation,
and gives the basis of all models mentioned above.
He also performed the ML estimation in the Weibull CE model with
the inverse power law on multiple-step SSALT data
concerning time to breakdown of an electrical insulation.
Hirose
(1996)
has proposed a generalized Weibull CE model,
which has a threshold parameter.

As we have seen, there are many kinds of studies about SSALTs
on the basis of the CE model and these provide significant
understanding of the Weibull and exponential models and SSALTs.
However, the validity of the models is not necessarily
clear
(Bagdanavicius, 1978; Nelson, 1980, 1990).
In this article we devote ourselves to considering the
following questions:
\begin{itemize}
\item
Can the models really express degradation of products?
\item
If so, what condition on the parameters is necessary for it?
\item
When a random variable obeys the Weibull CE model with a threshold parameter, 
how do its average and standard deviation behave under a condition?
\end{itemize}

In Section \ref{sec:ex_miner} we introduce the CE model with a threshold
parameter, which is a generalization of the CE model provided by Nelson.
After giving the Weibull CE model
with a threshold parameter
in Section \ref{sec:model}, we analyze it in Section \ref{sec:properties}.
In Section \ref{sec:tech} we give an ML estimation procedure.
In Section \ref{sec:example} we illustrate an example of the ML
estimation and a goodness of fit test
on an actual data set
and lastly give the conclusions.


%
\section{Cumulative exposure model}
\label{sec:ex_miner}
We construct a generalized CE model with the help of the CE model proposed
by Nelson
(1980, 1990),
whose model gives the 
distribution function of a random variable
for failure time. Although the general model is obtained
in a similar way to Nelson's, it differs in having a
threshold parameter that decides whether a specimen is
influenced by stress or not.

The assumptions to obtain the CE model were given
by Nelson as follows:
\begin{enumerate}
\renewcommand{\labelenumi}{\roman{enumi})}
\item
The remaining life of specimens depends only on the current
cumulative fraction accumulated.
\item
If held at the current stress, survivors will fail according
to the distribution function for that stress but starting
at the previously accumulated fraction failed.
\end{enumerate}

Using a distribution function $F$ of a non-negative random
variable with an explanatory
variable $V$ and a threshold $V_{th}$,
we construct the distribution function $G$
of a random variable $T$ for failure time
in a sequential way.
Denote by $V_{i}$ a stress that a specimen is subjected to in
an interval $(t_{i-1},t_{i}]$ (i=1,2,\ldots).

First of all, we define $G$ by
\begin{equation}
G(t)\defeq
\left\{
\begin{array}{ll}
F(t-t_{0};V_{1}) & (V_{1}>V_{th}),\\
0 & (V_{1}\leq V_{th})
\end{array}
\right.
\label{eq:Gt1}
\end{equation}
for $t_{0}\leq t\leq t_{1}$.

Next, for $t_{1}< t\leq t_{2}$ we define
\begin{equation}
G(t)\defeq
\left\{
\begin{array}{lr}
 F(t-t_{1}+s_{1};V_{2}) & (V_{2}>V_{th}),\\
F(s_{1};V_{2}) & (V_{2}\leq V_{th}).
\end{array}
\right.
\label{eq:Gt2}
\end{equation}
Here, according to Assumption ii), $s_{1}$ is a positive
value satisfying $G(t_{1})=F(s_{1};V_{2})$.

Similarly, for $t_{i-1}<t\leq t_{i}$ we define
\[
G(t)\defeq
\left\{
\begin{array}{lr}
F(t-t_{i-1}+s_{i-1};V_{i}) & (V_{i}>V_{th}),\\
F(s_{i-1};V_{i}) & (V_{i}\leq V_{th}),
\end{array}
\right.
\]
where $s_{i-1}$ is a positive value satisfying
$G(t_{i-1})=F(s_{i-1};V_{i})$.

Also in
Hirose (1996),
a similar formulation was given, provided
that the step stress at the present time is not lower than
that at the past time, which means $V_{m}\leq V_{th}$ holds
for any $m<i$ when $V_{i}\leq V_{th}$.
Actually, in his formulation
\[
G(t)\defeq
\left\{
\begin{array}{lr}
F(t-t_{i-1}+s_{i-1};V_{i}) & (V_{i}>V_{th}),\\
0 & (V_{i}\leq V_{th})
\end{array}
\right.
\]
for $t_{i-1}<t\leq t_{i}$.
Note that our formulation is more
general.


%
\section{Model for SSALT}
\label{sec:model}
We deal with a multiple-step SSALT
under the condition that specimens were subjected
to a normal level of stress and did not fail before the test.
As we will see in the next section, this setting
has a possibility of 
throwing light on new aspects concerning
the Weibull and exponential CE models.
In this section, first we introduce the multiple-step
SSALT and second we give the Weibull CE model under the condition.

\subsection{Multiple-step SSALT}
\label{subsec:step-test}
During the multiple-step SSALT, specimens are subjected
to successively higher levels of stress as follows.
After a specimen was used at a normal level
of stress, it is subjected  to an initial level of stress
for a predetermined time interval at the first stage in the test.
If it does not fail, it is subjected to a higher level of stress
for a predetermined time interval at the next stage.
In analogy, it is repeatedly subjected to higher levels of stress
until it fails. The other specimens are tested similarly.
The pattern of stress levels and time intervals is the same
for all specimens.

\subsection{Weibull CE Model with a threshold parameter}
We construct the Weibull CE model by combining the basic CE model
in Section \ref{sec:ex_miner} and
a Weibull law, in which
the underlying distribution is the two parameter Weibull and
the scale parameter is replaced with the function of an explanatory
variable. 
In addition, we give the log likelihood function
in the case when the step-stress
data are given under the condition mentioned above.

Let $V$ be an explanatory variable and a function of it
$\phi(V)$. When the scale parameter is replaced with
$\phi(V)$ in the Weibull distribution function, the distribution
is given by
\[
W(t;V) = 1-\exp\left[-\left(\frac{t}{\phi(V)}\right)^{\beta}\right].
\]
We use this as $F$ in Section \ref{sec:ex_miner}.

Denote by $V_{s}$ and $T_{s}$
a normal level of stress and
the length of the time interval during that
a specimen is used before the test, respectively.
In addition, denote by $V_{i}$ the stress that a
specimen is subjected to at the $(i-1)$-st stage
in the test, and let be $t_{i-1}$ the start time of the
stage $(i=2,3,4,\ldots)$. Since the level of stress
becomes higher as the stage in the test advances,
the relationship $V_{i} < V_{j}$ holds
when $2\leq i<j$.

Before we consider the Weibull CE model under the condition
mentioned in the first two lines of Section \ref{sec:model},
as preliminaries, let us consider the model without assuming
the condition holds. When we denote by $t_{1}$ the start time of the
test and set $t_{0}$ at $t_{1}-T_{s}$,
(\ref{eq:Gt1}) gives
\[
G(t_{1})=
\left\{
\begin{array}{ll}
\dstyle{1-\exp\left[-\left(\frac{T_{s}}{\phi(V_{s})}\right)^{\beta}\right]}
&\makebox[2em]{} (V_{s}>V_{th}),\\
0 &\makebox[2em]{} (V_{s}\leq V_{th}).
\end{array}
\right.
\]
This and
\[
W(s_{1};V_{2})
= 1-\exp\left[-\left(\frac{s_{1}}{\phi(V_{2})}\right)^{\beta}\right]
\]
yield
\[
s_{1}=
\left\{
\begin{array}{ll}
\dstyle{\frac{T_{s}}{\phi(V_{s})}\phi(V_{2})}
&\makebox[2em]{} (V_{s}>V_{th}),\\
0 &\makebox[2em]{} (V_{s}\leq V_{th}).
\end{array}
\right.
\]
Hence, (\ref{eq:Gt2}) gives
\[
G(t_{2})=
\left\{
\begin{array}{ll}
\dstyle{
  1-\exp\left[-\left(
      \frac{T_{s}}{\phi(V_{s})}+\frac{t_{2}-t_{1}}{\phi(V_{2})}
    \right)^{\beta}\right]
  \raisebox{0pt}[0pt][20pt]{}
}
&\makebox[2em]{} (V_{s}>V_{th},\ V_{2}>V_{th}),\\
\dstyle{
  1-\exp\left[-\left(
      \frac{T_{s}}{\phi(V_{s})}
    \right)^{\beta}\right]
}
&\makebox[2em]{} (V_{s}>V_{th},\ V_{2}\leq V_{th}),\\
\dstyle{
  1-\exp\left[-\left(
      \frac{t_{2}-t_{1}}{\phi(V_{2})}
    \right)^{\beta}\right]
}
&\makebox[2em]{} (V_{s}\leq V_{th},\ V_{2}> V_{th}),\\
0 &\makebox[2em]{} (V_{s}\leq V_{th},\ V_{2}\leq V_{th}).
\end{array}
\right.
\]
By repeating similar calculations, we can obtain
the cumulative distribution function 
$G$ in the $(i-1)$-st stage:
\begin{equation}
G(t) = 1- \exp
\left[-\varepsilon^{\beta}(t)\right],\qquad t_{i-1}< t\leq t_{i},
\label{eq:distG}
\end{equation}
where
\begin{equation}
\varepsilon(t)  \defeq
\left\{
\begin{array}{ll}
\dstyle{
  \frac{T_{s}}{\phi(V_{s})}
  +\frac{t_{k}-t_{k-1}}{\phi(V_{k})}
  + \cdots +\frac{t_{i-1}-t_{i-2}}{\phi(V_{i-1})}
  +\frac{t-t_{i-1}}{\phi(V_{i})}
}
&\makebox[2em]{} (V_{s}>V_{th})
\raisebox{0pt}[0pt][20pt]{},\\
\dstyle{
  \frac{t_{k}-t_{k-1}}{\phi(V_{k})}
  + \cdots +\frac{t_{i-1}-t_{i-2}}{\phi(V_{i-1})}
  +\frac{t-t_{i-1}}{\phi(V_{i})}
}
&\makebox[2em]{} (V_{s}\leq V_{th})
\end{array}
\right.
\label{eq:Ep}
\end{equation}
when
\[
V_{2}<V_{3}<\cdots<V_{k-1}\leq V_{th}< V_{k}
<\cdots<V_{i-1}< V_{i}.
\]
Since we are interested in the case of degradation
of products under a normal level of stress, 
we assume $V_{s}>V_{th}$ in the sequel.

Next, we seek our target, that is, the cumulative distribution
function under the condition that a specimen
was subjected to a normal level of stress and
did not fail before the test.
From the statements above, the function is as follows:
\begin{eqnarray}
G(t|t> t_{1})&=&
\frac{
\bigl\{G(t)-G(t_{1})\bigr\}I_{\{t> t_{1}\}}(t)
}{1-G(t_{1})}
\nonumber\\
&=&
\left\{
1-\exp\left[\varepsilon^{\beta}(t_{1})-\varepsilon^{\beta}(t)\right]
\right\}I_{\{t> t_{1}\}}(t),
\label{eq:condDist}
\end{eqnarray}
where
\[
I_{\{t> t_{1}\}}(t)\defeq
\left\{
\begin{array}{cc}
1 & (t> t_{1}), \\
0 & (t\leq t_{1}).
\end{array}
\right.
\]
Thus, the log likelihood function $\ln L$ under the condition
is expressed by the following:
if we denote by $N$ and $l$ the sample size and
the level of the stage at which a specimen fails,
and use superscript $(j)$ to show that a variable
is related to the $j$-th specimen,
\begin{equation}
\ln L = \sum_{j=1}^{N}\ln
\left\{
\exp\left(-\varepsilon^{\beta}(t_{l-1}^{(j)};T_{s}^{(j)})\right)
-\exp\left(-\varepsilon^{\beta}(t_{l}^{(j)};T_{s}^{(j)})\right)
\right\}
+\sum_{j=1}^{N}
\varepsilon^{\beta}(t_{1}^{(j)};T_{s}^{(j)}),
\label{eq:likeli}
\end{equation}
where we express $\varepsilon(t)$ by $\varepsilon(t;T_{s})$ in order
to show clearly that each specimen has each $T_{s}$.
%


%
%
\section{Statistical properties}
\label{sec:properties}
We consider the statistical properties of the model under the condition
mentioned in the previous section.
First we state the role of the shape parameter
$\beta$ in the distribution function (\ref{eq:condDist})
and second we investigate the relationship between statistical quantities
and the values of parameters after we
simplify the model without loss of generality.
In the sequel we express $G(t|t>t_{1})$
by $G(t|t>t_{1};T_{s})$ when it is necessary to show clearly the length
of the time interval during that a specimen is used
before the test.
In analogy, we express $G(t|t>t_{1})$
by $G(t|t>t_{1};V_{s})$ when it is necessary to show clearly the normal
stress that a specimen is subjected to
before the test. Since we are interested in
elapsed time from the start time of test,
in the sequel we suppose that $t_{1}$ is
the base point in time. That is,
we may consider $t_{1}$ equal to $0$.

Depending on the magnitude of $\beta$, the distribution function
has a different aspect as follows.

\begin{Lem_eng}
\label{lem:gfuncTS}
Assume that $T_{a}<T_{b}$.
Then, the following holds for $t>t_{k-1}$
such that $V_{k-1}\leq V_{th}<V_{k}$.
\begin{enumerate}
\renewcommand{\labelenumi}{\roman{enumi})}
\item
If $0<\beta<1$,
\[
G(t|t>t_{1};T_{a})>G(t|t>t_{1};T_{b}).
\]
\item
If $\beta=1$,
\[
G(t|t>t_{1};T_{a})=G(t|t>t_{1};T_{b}).
\]
\item
If $\beta>1$,
\[
G(t|t>t_{1};T_{a})<G(t|t>t_{1};T_{b}).
\]
\hfill \bsquare
\end{enumerate}
\end{Lem_eng}

\noindent {\it Proof.}
The substitutions of $\varepsilon(t_{1})=T_{s}/\phi(V_{s})$
and (\ref{eq:Ep}) into (\ref{eq:condDist}) yield
\[
G(t|t>t_{1};T_{s})=1-\exp\left[
\left\{
\frac{T_{s}}{\phi(V_{s})}
\right\}^{\beta}
-
\left\{
\frac{T_{s}}{\phi(V_{s})}+
\sum_{m=k}^{i-1}\frac{t_{m}-t_{m-1}}{\phi(V_{m})}
+\frac{t-t_{i-1}}{\phi(V_{i})}
\right\}^{\beta}
\right].
\]
By differentiating this with respect to $T_{s}$ and arranging it, we find
\[
\frac{\partial G(t|t>t_{1};T_{s})}{\partial T_{s}}
=
\frac{\beta}{\phi(V_{s})}
\left[
\varepsilon^{\beta-1}(t)
-\varepsilon^{\beta-1}(t_{1})
\right]
\exp
\left[
\varepsilon^{\beta}(t_{1})
-\varepsilon^{\beta}(t)
\right].
\]
Noting $\varepsilon(t)>\varepsilon(t_{1})$, we can see
\begin{enumerate}
\renewcommand{\labelenumi}{\roman{enumi})}
\item
if $0<\beta<1$,\\
$G(t|t>t_{1};T_{s})$ is a strictly decreasing function
of $T_{s}$ since $\partial G(t|t>t_{1};T_{s})/\partial T_{s}<0$,
\item
if $\beta=1$,\\
$G(t|t>t_{1};T_{s})$ does not depend on $T_{s}$
since $\partial G(t|t>t_{1};T_{s})/\partial T_{s}=0$,
\item
if $\beta>1$,\\
$G(t|t>t_{1};T_{s})$ is a strictly increasing function
of $T_{s}$ since $\partial G(t|t>t_{1};T_{s})/\partial T_{s}>0$.
\end{enumerate}
This completes the proof.
\hfill $\Box$

The statement i) in the lemma means that specimens become  more
durable as they are used longer before the test.
This is clearly irrational. Thus, in this sense any value in
$(0,1)$ is inadmissible for $\beta$. The statement ii) deals with
a situation when the underlying distribution is exponential.
It indicates that the CE model inherits the memoryless
property from the exponential distribution. The statement iii)
expresses the most realistic situation, in which the durability
of specimens decreases as the the duration of their use becomes
longer before the test.

In a similar fashion, we can obtain the following lemma.

\begin{Lem_eng}
\label{lem:gfuncVS}
Assume that $(V_{th}<)V_{a}<V_{b}$.
Then, the following holds for $t>t_{k-1}$
such that $V_{k-1}\leq V_{th}<V_{k}$.
\begin{enumerate}
\renewcommand{\labelenumi}{\roman{enumi})}
\item
If $0<\beta<1$,
\[
G(t|t>t_{1};V_{a})>G(t|t>t_{1};V_{b}).
\]
\item
If $\beta=1$,
\[
G(t|t>t_{1};V_{a})=G(t|t>t_{1};V_{b}).
\]
\item
If $\beta>1$,
\[
G(t|t>t_{1};V_{a})<G(t|t>t_{1};V_{b}).
\]
\hfill \bsquare
\end{enumerate}
\end{Lem_eng}

From this lemma, we can know a similar fact to Lemma \ref{lem:gfuncTS}.
Especially, note that the statement iii) expresses the most
realistic situation, in which the durability of specimens decreases
as the normal stress imposed before the test becomes higher.

In the sequel we assume the inverse power law in $\phi$ for
$V>V_{th}$:
\begin{equation}
\phi(V) = \frac{K}{(V-V_{th})^n},
\label{eq:phi}
\end{equation}
where $K$ and $n$ are positive parameters and $V_{th}$ is a non-negative
parameter.
In addition, we assume that the length of the time interval and
the breadth of upsurge of stress are constant in the test.
That is, we set
\[
\del{t}\defeq t_{i}-t_{i-1},\quad
\del{V}\defeq V_{i+1}-V_{i}\quad(i=2,3,\ldots)
\quad {\rm and}\quad V_{2}=\del{V}.
\]
Let us simplify (\ref{eq:Ep}) and seek the
expectation and second moment of a random variable
obeying (\ref{eq:condDist}).

By using the above constants and rewriting
(\ref{eq:Ep}) and (\ref{eq:phi}), we can obtain
\begin{eqnarray}
\varepsilon(t) &=&
\frac{\del{t}}{\phi(V_{s})}\tilde{T}_{s}
+\frac{\del{t}}{\phi(V_{k})}
+ \cdots +\frac{\del{t}}{\phi(V_{i-1})}
+\frac{\del{t}}{\phi(V_{i})}\frac{t-t_{i-1}}{\del{t}},
\label{eq:simEp}
\\
\frac{\del{t}}{\phi(V_{s})} &=& \frac{(1-\tilde{V}_{th})^n}{\tilde{K}},
\quad
\frac{\del{t}}{\phi(V_{m})} =
\frac{((m-1)\del{\tilde{V}}-\tilde{V}_{th})^n}{\tilde{K}}
\quad (m=k,k+1,\ldots,i),
\qquad\quad
\label{eq:simPhi}
\end{eqnarray}
where
\[
\tilde{T}_{s}\defeq\frac{T_{s}}{\del{t}},
\quad \tilde{K}\defeq\frac{K}{(\del{t}V_{s})^{n}},
\quad \del{\tilde{V}}\defeq\frac{\del{V}}{V_{s}},
\quad \tilde{V}_{th}\defeq\frac{V_{th}}{V_{s}}
\]
and $(k-2)\del{\tilde{V}}\leq\tilde{V}_{th}<(k-1)\del{\tilde{V}}$
holds.
Here, remark that $\tilde{K}$ and $\tilde{V}_{th}$ are parameters
to be estimated while $\tilde{T}_{s}$ and $\del{\tilde{V}}$ are
quantities to be prespecified in order to decide a concrete model.
These expressions indicate that we can take
$\del{t}$ and $V_{s}$ as a unit of time and a unit of stress,
respectively. Besides, we can suppose that $0\leq \tilde{V}_{th}<1$
when we deal with the case that $V_{s}>V_{th}$.

By means of a similar procedure and the arrangement of expressions, we
can obtain another $\varepsilon(t)$ for a different
stress $V^{\prime}_{s}$, say $\varepsilon^{\prime}(t)$,
in the following form:
\[
\varepsilon^{\prime}(t) =
\frac{\del{t}}{\phi(V^{\prime}_{s})}\tilde{T}_{s}
+\frac{\del{t}}{\phi(V_{k})}
+ \cdots +\frac{\del{t}}{\phi(V_{i-1})}
+\frac{\del{t}}{\phi(V_{i})}\frac{t-t_{i-1}}{\del{t}},
\]
where
\[
\frac{\del{t}}{\phi(V^{\prime}_{s})}
= \frac{((V^{\prime}_{s}/V_{s})-\tilde{V}_{th})^n}{\tilde{K}}.
\]
Note that only the first terms in the right-hand sides
differ in
the expressions of $\varepsilon(t)$ and $\varepsilon^{\prime}(t)$.
Thus, once we obtain the values of parameters,
we can decide the distribution function in the
case of another stress $V^{\prime}_{s}$ ($>V_{s}$)
by replacing only the first term in the right-hand side of (\ref{eq:simEp}).

Let us seek the expectation of a random variable $T$ obeying
(\ref{eq:condDist}). We first seek the following
conditional expectation as preliminaries: for $m>k$,
\begin{eqnarray*}
E[T|T\leq t_{m}]
&=&\sum_{i=k}^{m}\int_{t_{i-1}}^{t_{i}}
t\frac{\partial}{\partial t}G(t|t>t_{1})\d t
\\
&=&\sum_{i=k}^{m}\left[
\Bigl[
tG(t|t>t_{1})
\Bigr]_{t_{i-1}}^{t_{i}}
-\int_{t_{i-1}}^{t_{i}}G(t|t>t_{1})\d t
\right]
\\
&=&-t_{m}\exp[\varepsilon^{\beta}(t_{1})-\varepsilon^{\beta}(t_{m})]
+t_{k-1}
+\sum_{i=k}^{m}\int_{t_{i-1}}^{t_{i}}
\exp[\varepsilon^{\beta}(t_{1})-\varepsilon^{\beta}(t)]\d t.
\end{eqnarray*}
In the last line of this equation the relationship
$\varepsilon(t_{1})=\varepsilon(t_{k-1})$ is used,
which holds by (\ref{eq:Ep}).

When we denote by $m_{0}$ a positive integer
such that $1/\phi(V_{i})<1$ holds for any $i>m_{0}$, we can see that
\begin{eqnarray*}
t_{m}\exp[-\varepsilon^{\beta}(t_{m})]
&=&
(m-1)\del{t}\exp\left[
-\left\{\frac{\del{t}}{\phi(V_{s})}\tilde{T}_{s}
+\sum_{i=k}^{m}\frac{\del{t}}{\phi(V_{i})}\right\}^{\beta}
\right]
\\
&<&(m-1)\del{t}\exp\left[
-\left\{\frac{\del{t}}{\phi(V_{s})}\tilde{T}_{s}
+\sum_{i=k}^{m_{0}}\frac{\del{t}}{\phi(V_{i})}+(m-m_{0})\del{t}
\right\}^{\beta}
\right]
\end{eqnarray*}
and the right-hand side converges to $0$ as $m\to\infty$.

From the things above and
$\dstyle{E[T]=\lim_{m\to\infty}E[T|T\leq t_{m}]}$, we obtain
\begin{eqnarray}
\frac{E[T]}{\del{t}}
&=&\frac{t_{k-1}}{\del{t}}
+\frac{1}{\del{t}}\sum_{i=k}^{\infty}\int_{t_{i-1}}^{t_{i}}
\exp[\varepsilon^{\beta}(t_{1})-\varepsilon^{\beta}(t)]\d t
\nonumber\\
&=&\frac{t_{k-1}}{\del{t}}
+\frac{1}{\beta}\sum_{i=k}^{\infty}\frac{\phi(V_{i})}{\del{t}}
\Bigl\{-A(t_{i})+A(t_{i-1})\Bigr\}
\label{eq:mean}
\end{eqnarray}
as the expectation in the case that $\del{t}$ is used as a unit of time.
Here,
\[
A(t)\defeq
\exp[\varepsilon^{\beta}(t_{1})-\varepsilon^{\beta}(t)]
\int_{0}^{\infty}
\left\{u+\varepsilon^{\beta}(t)\right\}^{1/\beta-1}
\e^{-u}\d u.
\]
The expression in the right-hand side
of (\ref{eq:mean}) is useful for stable numerical calculations
when $T_{s}$ takes a large value.

In a similar fashion we obtain
\begin{eqnarray}
\frac{E[T^{2}]}{(\del{t})^{2}}
&=&\left(\frac{t_{k-1}}{\del{t}}\right)^{2}
+\frac{2}{(\del{t})^{2}}\sum_{i=k}^{\infty}\int_{t_{i-1}}^{t_{i}}
t\exp[\varepsilon^{\beta}(t_{1})-\varepsilon^{\beta}(t)]\d t
\nonumber\\
&=&\left(\frac{t_{k-1}}{\del{t}}\right)^{2}
+\frac{2}{\beta}\sum_{i=k}^{\infty}\left(\frac{\phi(V_{i})}{\del{t}}\right)^{2}
\Bigl\{-B_{i}(t_{i})+B_{i}(t_{i-1})\Bigr\}
\label{eq:moment}
\end{eqnarray}
as the second moment in the case that $\del{t}$ is used as a unit of time.
Here,
\begin{eqnarray*}
B_{i}(t)&\defeq&
\exp[\varepsilon^{\beta}(t_{1})-\varepsilon^{\beta}(t)]
\\
&&
\makebox[2em]{}
\times
\int_{0}^{\infty}
\left\{\Bigl(u+\varepsilon^{\beta}(t)\Bigr)^{1/\beta}
-\left(\varepsilon(t_{i-1})-\frac{t_{i-1}-t_{1}}{\phi(V_{i})}\right)
\right\}
\left(u+\varepsilon^{\beta}(t)\right)^{1/\beta-1}
\e^{-u}\d u.
\end{eqnarray*}
Using (\ref{eq:mean}) and (\ref{eq:moment}), we can calculate the
mean and standard deviation of $T/\del{t}$ for the 
parameter values in Table 1.
The results are shown on Figs 1 and 2.
In these figures
we show the difference in the pair of the values of
$\tilde{V}_{th}$ and $n$ by the combination of the sort of line
and the thickness or the color of line. That is,
the solid, dash or dotted line
means that $\tilde{V}_{th}=0$, $0.5$ or
$0.9$, respectively.
On the other hand, the thick, normal or gray one
means $n=1$, $2$ or $3$, respectively.
%
%
\begin{table}[t]
\renewcommand{\thetable}{\arabic{table}}
\caption{Parameter values}
\label{tab:Para_values}
\begin{center}
\begin{tabular}{c|c|c|c|c}
$\tilde{K}$ & $\del{\tilde{V}}$ & $\tilde{V}_{th}$ & $\beta$ & $n$ \\
\hline
$10^{3}$, $10^{4}$, $10^{5}$ & $0.39$
& $0$, $0.5$, $0.9$ & $2$, $3$ & $1$, $2$, $3$
\end{tabular}
\end{center}
\end{table}

\clearpage
%
%
\begin{figure}[h]
\begin{center}
\unitlength 1cm
\begin{picture}(16,18.5)(0,0)
\put(0.5,0.5){\includegraphics[width=7cm]{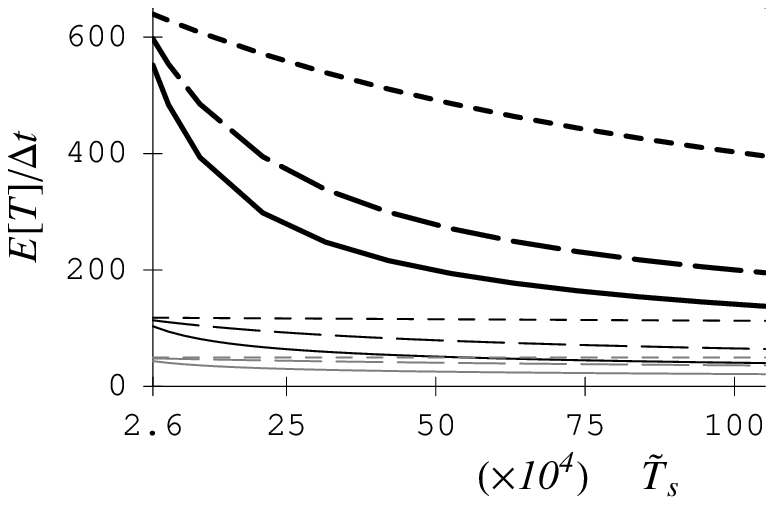}}
\put(8.5,0.5){\includegraphics[width=7cm]{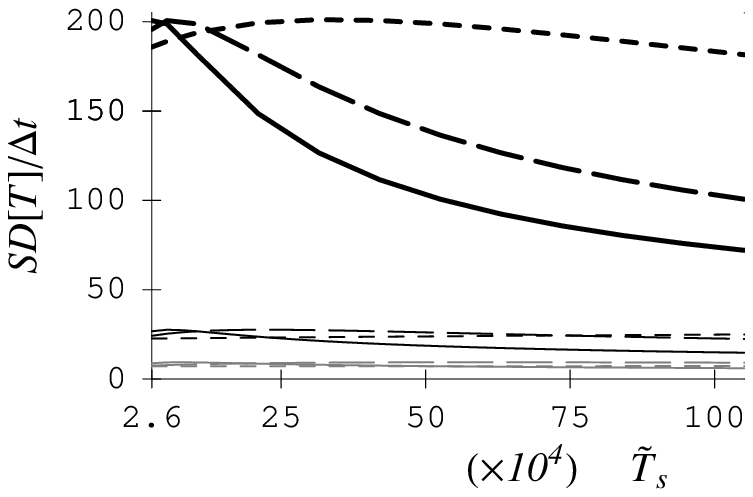}}
\put(0.5,7){\includegraphics[width=7cm]{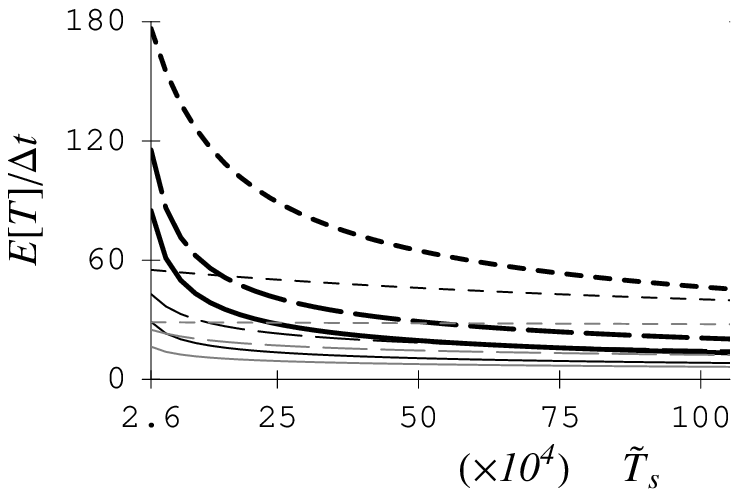}}
\put(8.5,7){\includegraphics[width=7cm]{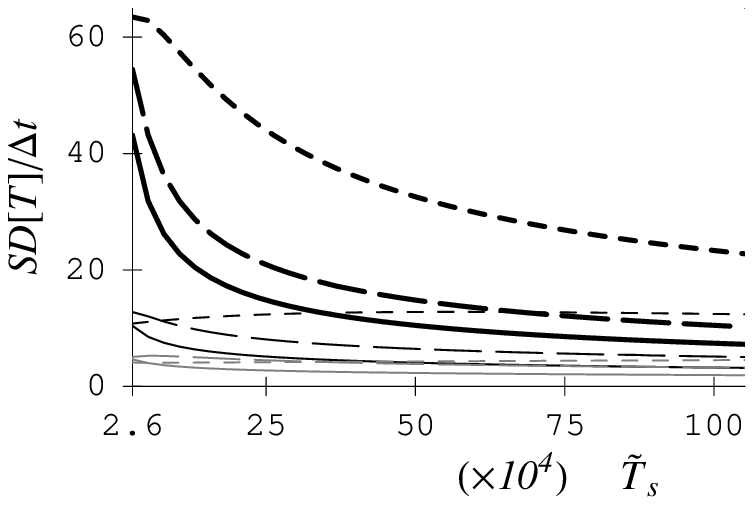}}
\put(0.5,13.5){\includegraphics[width=7cm]{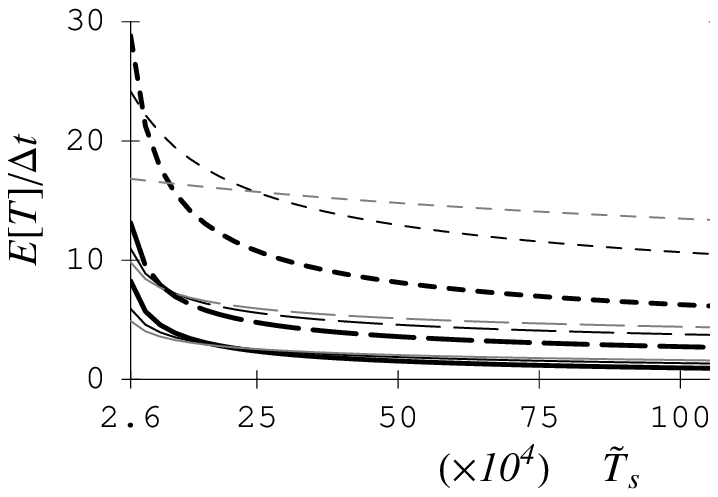}}
\put(8.5,13.5){\includegraphics[width=7cm]{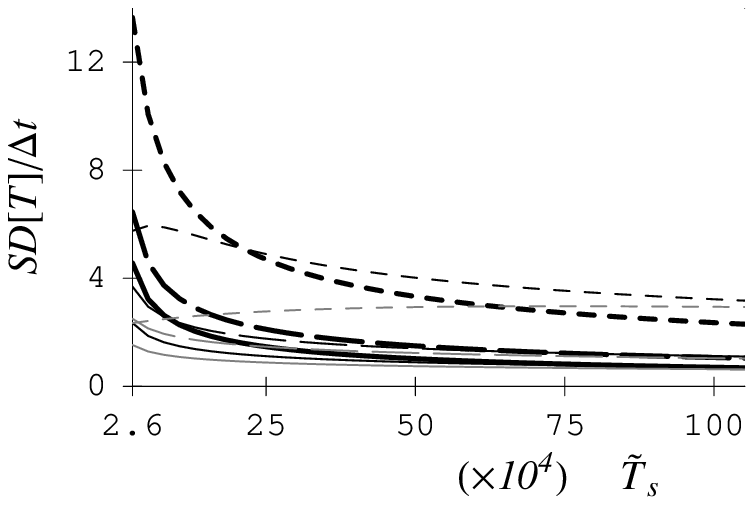}}
\put(6.8,12.9){The case of $\tilde{K}=10^{3}$}
\put(6.8,6.4){The case of $\tilde{K}=10^{4}$}
\put(6.8,-0.1){The case of $\tilde{K}=10^{5}$}
\end{picture}
\end{center}
\caption{The mean and standard deviation of $T/\del{t}$ when $\beta=2$}
\label{fig:b2}
\end{figure}
%

%
\clearpage
%
%
\begin{figure}[t]
\begin{center}
\unitlength 1cm
\begin{picture}(16,18.5)(0,0)
\put(0.5,0.5){\includegraphics[width=7cm]{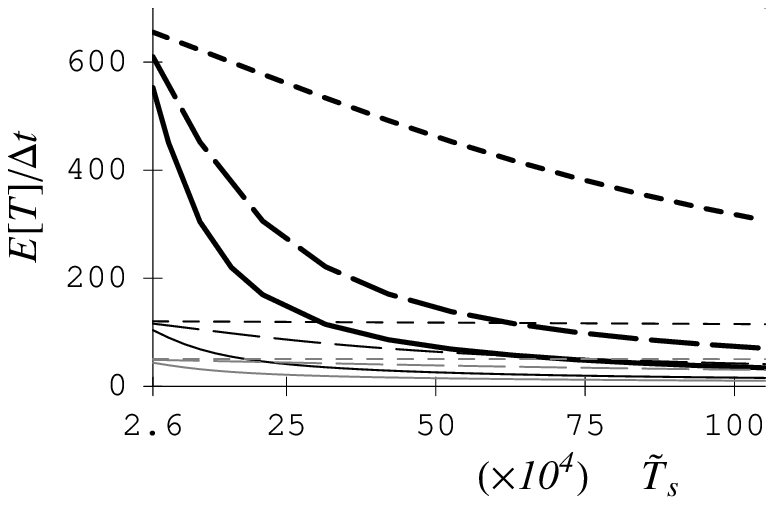}}
\put(8.5,0.5){\includegraphics[width=7cm]{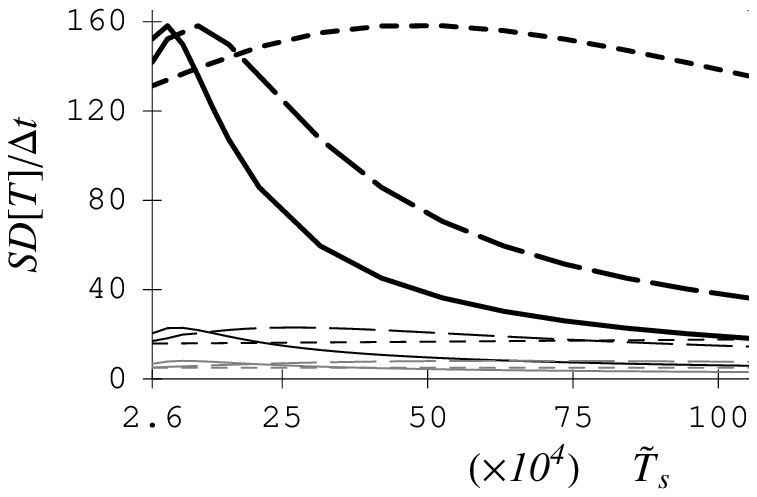}}
\put(0.5,7){\includegraphics[width=7cm]{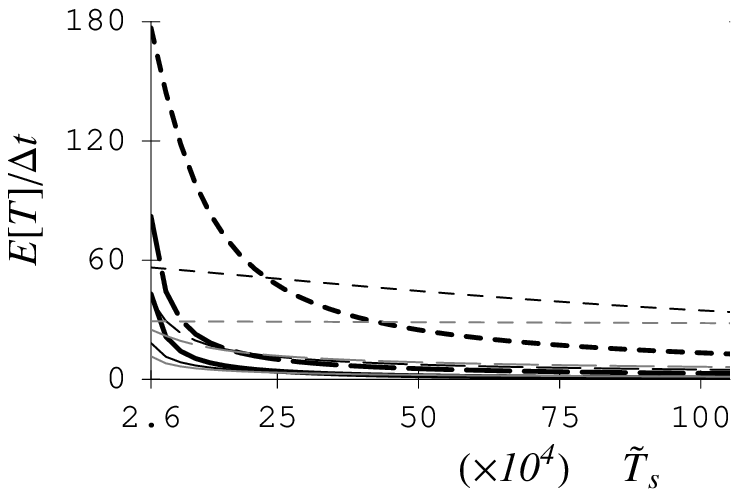}}
\put(8.5,7){\includegraphics[width=7cm]{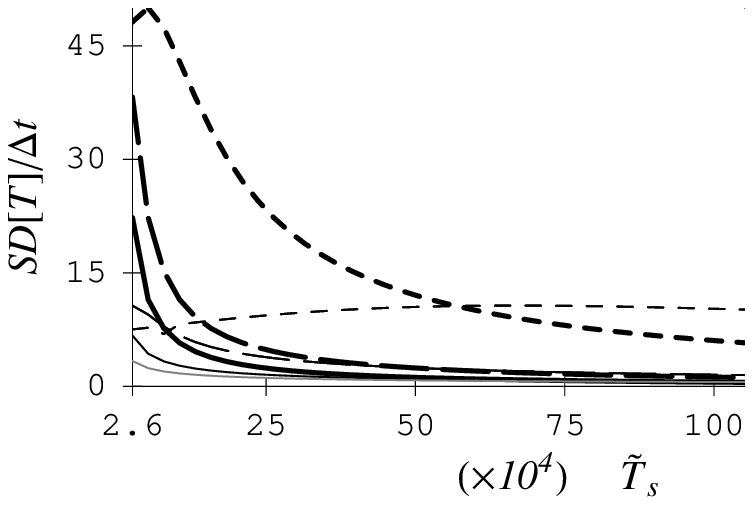}}
\put(0.5,13.5){\includegraphics[width=7cm]{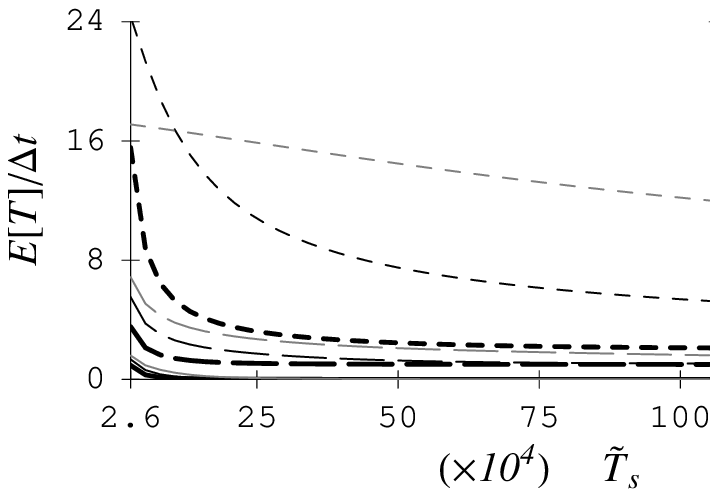}}
\put(8.5,13.5){\includegraphics[width=7cm]{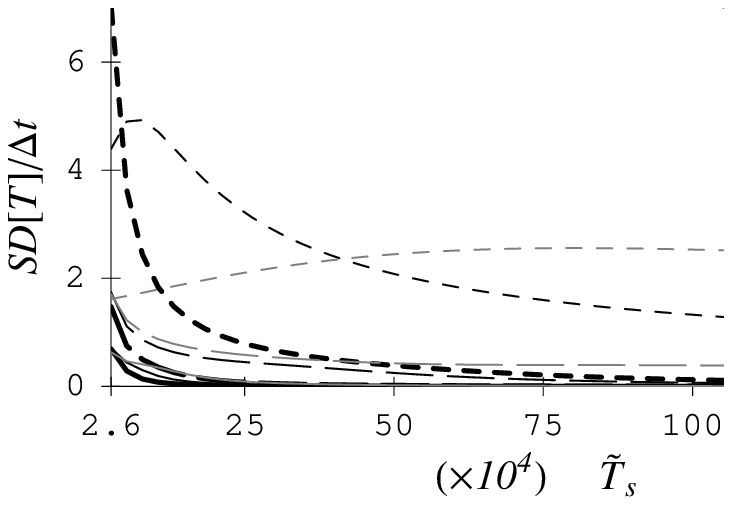}}
\put(6.8,12.9){The case of $\tilde{K}=10^{3}$}
\put(6.8,6.4){The case of $\tilde{K}=10^{4}$}
\put(6.8,-0.1){The case of $\tilde{K}=10^{5}$}
\end{picture}
\end{center}
\caption{The mean and standard deviation of $T/\del{t}$ when $\beta=3$}
\label{fig:b3}
\end{figure}

\clearpage
%

%

%
\section{Estimation procedure}
\label{sec:tech}

We state the way of seeking the ML estimates of the parameters
in (\ref{eq:condDist}), (\ref{eq:simEp}) and (\ref{eq:simPhi}).
Some techniques below help us to obtain the estimates
numerically stably.

We use 
a new parameter $\zeta$ defined by $\tilde{K}/\tilde{K}_{0}$
for a constant $\tilde{K}_{0}$
instead of $\tilde{K}$.
The reason is because
only the parameter
$\tilde{K}$ possibly has an estimate that is much larger
than those of the other parameters.

By differentiating (\ref{eq:likeli}) with respect to each parameter
and arranging each equation, we can obtain the likelihood equations
in the following simplified form:
\begin{equation}
\sum_{j=1}^{N}\frac{\lambda_{\theta}^{(j)}}{d_{j}}
+\sum_{j=1}^{N}\delta_{\theta}(t_{1}^{(j)};\tilde{T}_{s}^{(j)})=0,
\quad \theta\in\{\beta,n,\zeta,\tilde{V}_{th}\},
\label{eq:likeliEQ}
\end{equation}
where
\begin{eqnarray*}
d_{j}&\defeq&
\exp\left(-\varepsilon^{\beta}(t_{l-1}^{(j)};\tilde{T}_{s}^{(j)})\right)
-\exp\left(-\varepsilon^{\beta}(t_{l}^{(j)};\tilde{T}_{s}^{(j)})\right),
\\
\lambda_{\theta}^{(j)}&\defeq&
-\delta_{\theta}(t_{l-1}^{(j)};\tilde{T}_{s}^{(j)})
\exp\left(-\varepsilon^{\beta}(t_{l-1}^{(j)};\tilde{T}_{s}^{(j)})\right)
+\delta_{\theta}(t_{l}^{(j)};\tilde{T}_{s}^{(j)})
\exp\left(-\varepsilon^{\beta}(t_{l}^{(j)};\tilde{T}_{s}^{(j)})\right),
\\
\delta_{\beta}(t_{i};\tilde{T}_{s})&\defeq&
\varepsilon^{\beta}(t_{i};\tilde{T}_{s})\ln \varepsilon(t_{i};\tilde{T}_{s}),
\\
\delta_{n}(t_{i};\tilde{T}_{s})&\defeq&
\frac{1}{\tilde{K}_{0}}\varepsilon^{\beta-1}(t_{i};\tilde{T}_{s})
\left\{
\sum_{m=k}^{i}((m-1)\del{\tilde{V}}-\tilde{V}_{th})^{n}
\ln ((m-1)\del{\tilde{V}}-\tilde{V}_{th})
\right\},
\\
\delta_{\zeta}(t_{i};\tilde{T}_{s})&\defeq&
\varepsilon^{\beta}(t_{i};\tilde{T}_{s}),
\quad
\delta_{\tilde{V}_{th}}(t_{i};\tilde{T}_{s})\defeq
\varepsilon^{\beta-1}(t_{i};\tilde{T}_{s})
\frac{1}{\tilde{K}_{0}}
\sum_{m=k}^{i}((m-1)\del{\tilde{V}}-\tilde{V}_{th})^{n-1}.
\end{eqnarray*}
In the expressions above, note that 
$\varepsilon(t)$ in (\ref{eq:simEp}) is expressed
by $\varepsilon(t;\tilde{T}_{s})$ as usual.
We seek the zero of (\ref{eq:likeliEQ}) by means of the damped Newton method
(Bank and Rose, 1981) as follows.

In this model, 
the calculation for ML estimates 
is so sensitive that, depending on a vector of initial guesses, 
a sequence of approximate vectors by the damped Newton iteration
can converge to
a vector of estimates in which the estimate of $\beta$ is less 1
even if the value of the likelihood function is
not a maximum value. According to our observation, when
this phenomenon occurs, there is often
a tendency that the estimate of $\tilde{V}_{th}$
tends to 0. Thus, we adopt the strategy below.
\begin{enumerate}
  \item
    We seek the profile
    of (\ref{eq:likeli}) with respect to $\tilde{V}_{th}$. That is,
    while changing the value of $\tilde{V}_{th}$ from
    a value $\alpha_{0}$ to another value $\alpha_{1}$ in incremental steps,
    we seek the estimates of the other parameters in each step.
  \item
    Among the points on the part of the profile, we select the point at that
    the profile achieves its maximum, and then seek the ML estimates
    of all parameters simultaneously by using the point as a vector
    of initial values
    and performing the damped Newton iteration
    without fixing $\tilde{V}_{th}$.
\end{enumerate}
The derivatives of the expressions in the left-hand side of
(\ref{eq:likeliEQ}) are given in Appendix.


%
\section{Example}
\label{sec:example}

On the basis of the results obtained in the previous two sections,
we show an example of the ML estimation
in (\ref{eq:condDist}) on an actual data set. The data set
is part of the step-stress data on time to breakdown of cross-linked
polyethylene-insulated cables in
Hirose (1997).
We chose only data whose insulation class is 22kV (in 3-phase).
The reason is because data whose insulation class is 33kV
in the literature indicate
that their durability is higher than those for 22kV.
Thus, we concluded that we can not mix both of them.
Also, note Lemma \ref{lem:gfuncVS}. Moreover, we did not
incorporate into our sample data set
the data coming from the cables that passed 26 years
because their values are abnormal, compared with the others, and
they push down the value of the likelihood function extremely.
Finally, we
perform a goodness of fit test by utilizing the Monte Carlo method.

\subsection{Data set}
We show
the data set in Table 2.
For each specimen $j$,
the first column indicates the length of the time
interval during that the specimen is used before
the test, and
the second column indicates elapsed time from the start
time of test by the start time of the stage on which
a specimen fails. The unit of time is ten minutes.
The third column indicates the number of data in each row,
which is denoted by $N_{d}$.
The fourth and fifth columns indicate the average and
standard deviation of data, respectively,
in the case that an outrageous datum is not taken into account.
%

%
\subsection{Maximum likelihood estimation}
On the data set, we can see that $\del{\tilde{V}}=0.39$
because that $V_{s}=22$kV (in 3-phase) $=22/\sqrt{3}$kV
(in single phase) and $\del{V}=5$kV (in single phase).
As a constant for $\tilde{K}_{0}$,
we set $\tilde{K}_{0}=10^{4}$ and used as a vector of initial guesses
$(\beta,n,\zeta,\tilde{V}_{th})=(2.0,2,1,0.5)$.
These were roughly guessed from the comparison
between Fig. 1 or 2 and
the averages or standard deviations in Table 2.
In addition, for seeking the profile we set $\alpha_{0}=0.5$,
$\alpha_{1}=0.999$ and the increment of $\tilde{V}_{th}$
at $0.001$ in each step.

Ultimately, we can obtain
the following ML estimates:
\begin{equation}
\beta=5.016812,\quad n=1.603875,
\quad \zeta=0.548237,\quad \tilde{V}_{th}=0.944054.
\label{eq:paraV}
\end{equation}
Then, the value of the likelihood function (\ref{eq:likeli})
is $-244.4626$. The mean and some statistical quantities of $T$
and the scatter plot of data
are given on Fig. 3. In the
figure the solid line indicates the mean and
the upper or lower dotted line
indicates the mean plus or minus the standard
deviation, respectively. Each dot indicates
$(t_{l-1}^{(j)}+t_{l}^{(j)})/2$ for each sample $j$.

\clearpage
%
%
\begin{table}[h]
\renewcommand{\thefootnote}{\fnsymbol{footnote}}
\caption{Step-stress data}
\label{tab:data}
\begin{center}
\begin{tabular}{c|l|c|l|r}
$\tilde{T}_{s}^{(j)}$ & \multicolumn{1}{c|}{$t_{l-1}^{(j)}/\del{t}$} & $N_{d}$
& \multicolumn{1}{c|}{Ave.} & \multicolumn{1}{c}{SD} \\
\hline
157680 & 54, 56, 59, 64 & 4 & 58.3 & 4.3\\
473040 & 16 & 1 & 16 & \multicolumn{1}{c}{$\ast$} \\
578160 & 16, 16, 18, 19, 22, 24, 46, 48, 50 & 9 & 28.8 & 14.7\\
630720 & 24, 24, 28, 30, 32 & 5 & 27.6 & 3.6\\
735840 & 19, 20, 23, 35, 39, 42 & 6 & 30.0 & 10.2\\
788400 & 12, 12, 12, 13, 14, 14, 14, 17, 17 & 9 & 13.9 & 2.0\\
840960 & 12, 14, 14, 23, 24, 74\footnotemark[2] & 6 & 17.4 & 5.6\\
893520 & 16, 18, 20, 20, 26, 35 & 6 & 22.5 & 7.0\\
946080 & 16, 16, 16, 17, 18, 18, 20, 22, 26 & 9 & 18.8 & 3.4\\
998640 & 8, 10, 11, 12, 12, 12, 13, 13, 14, 15, 22 & 11& 12.9 & 3.6\\
1051200& 11, 12, 12 & 3 & 11.7 & 0.6\\
1156320& 11, 12, 12, 13, 14, 14 & 6 & 12.7 & 1.2
\end{tabular}
\end{center}
\vspace*{-2mm}
\hspace*{15mm}
{\footnotesize The marks \dag\ and $\ast$ mean outrageous and incomputable, respectively.}
\end{table}

%
%
\begin{figure}[h]
\begin{center}
\unitlength 1cm
\begin{picture}(8,5.5)(0,0)
\put(0.5,0.5){\includegraphics[width=7cm]{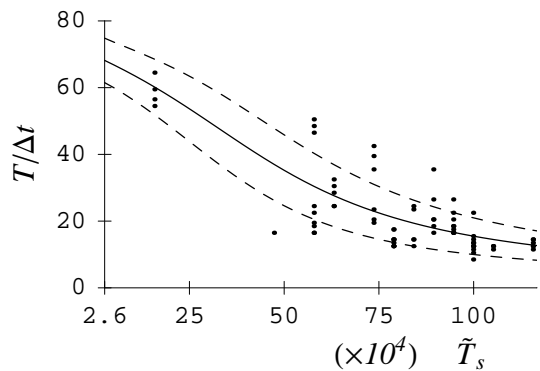}}
\end{picture}
\end{center}
\caption{Mean and scatter plot}
\label{fig:res}
\end{figure}
%

%
\subsection{Goodness of fit test}
After the values of the parameters were estimated and
a model was decided,
we are often concerned with testing its validity.
In the example, however,
we can not perform a chi-square goodness
of fit test on the grouped data because
the number of samples is too small
(Rao, 2002, p. 396).
Hence, in almost the same way as that Ross (2002, p. 206) adopted
in such a situation,
we test the hypothesis that
the model is consistent with the data set
when the parameters are of the values in (\ref{eq:paraV}).

\subsubsection*{Simulation conditions}
We performed Monte Carlo simulation under the simulation conditions
below.
\begin{itemize}
\item Setting for generating simulated data\\
We used as prespecified values
\[
\del{\tilde{V}}=0.39,\quad \tilde{K}_{0}=10^{4},\quad \tilde{T}_{s}
=157680,473040,\ldots,1156320,
\]
which are the same values as those for the data set in Table 2, and
as the true values of the parameters
\begin{equation}
\beta=5.016812,\quad n=1.603875,
\quad \zeta=0.548237,\quad \tilde{V}_{th}=0.944054,
\label{eq:trueVal}
\end{equation}
which come from (\ref{eq:paraV}).

~~~For each value of $\tilde{T}_{s}$, the number of simulated data
is the same as that in Table 2. The outrageous datum
in the table is, however, counted out. That is, for $\tilde{T}_{s}=840960$,
the number of simulated data is $5$ not $6$. Hence, the
total number of data is $74$ in a sample set.
\item Sample sets\\
Except sample sets where the ML estimates could not be obtained,
1000 sets of independent pseudo-random samples were considered.
\item Setting for estimation\\
On the stage to seek the profile we set $\alpha_{0}=0.85$,
$\alpha_{1}=0.999$ and the increment of $\tilde{V}_{th}$ at
0.001 in each step.
\end{itemize}

\subsubsection*{Procedure for generating random sets}
From (\ref{eq:Ep}), (\ref{eq:simEp}) and (\ref{eq:simPhi})
\[
\varepsilon(t_{1})=\frac{T_{s}}{\phi(V_{s})}
=\frac{(1-\tilde{V}_{th})^{n}}{\tilde{K}_{0}\zeta}\tilde{T}_{s}.
\]
When we set $u\defeq G(t|t>t_{1})$ for ($\infty>$) $t>t_{k-1}$,
this and (\ref{eq:condDist}) give
\begin{equation}
\varepsilon(t)=
\left[
-\ln(1-u)+\left(\frac{(1-\tilde{V}_{th})^{n}}{\tilde{K}_{0}\zeta}
\tilde{T}_{s}\right)^{\beta}
\right]^{1/\beta}.
\label{eq:q}
\end{equation}
Let us denote by $q$ the expression in the right-hand side.
Then, (\ref{eq:simEp}) and (\ref{eq:simPhi}) yield
\begin{equation}
\sum_{m=k}^{i-1}
\frac{((m-1)\del{\tilde{V}}-\tilde{V}_{th})^{n}}{\tilde{K}_{0}\zeta}
<q-\frac{(1-\tilde{V}_{th})^{n}}{\tilde{K}_{0}\zeta}\tilde{T}_{s}\leq
\sum_{m=k}^{i}
\frac{((m-1)\del{\tilde{V}}-\tilde{V}_{th})^{n}}{\tilde{K}_{0}\zeta}.
\label{ineq:q}
\end{equation}
Note that $u\in (0,1)$ and $i\geq k$ since $t_{k-1}<t<\infty$.
Consequently, the procedure for generating a pseudo-random set is
as follows:
\begin{enumerate}
\item generate a uniform random number and seek $q$ by (\ref{eq:q}),
\item find $i$ that satisfies (\ref{ineq:q}), where the $i$ is
  the stress level $l$ in which failure occurs,
\item calculate $t_{l-1}$,
\item repeat the three steps above 74 times.
\end{enumerate}
\subsubsection*{Procedure for a goodness of fit test}
In Table 2 we choose some
$\tilde{T}_{s}$'s and calculate the value of the
test statistic
\begin{equation}
T=\sum_{i=1}^{\kappa}\frac{(m_{i}-N_{d}p_{i})^{2}}{N_{d}p_{i}}
\label{eq:T}
\end{equation}
for each of the $\tilde{T}_{s}$'s.
Here, $m_{i}$ stands for the number of data in the $i$th subinterval
when the interval where failure time lies is divided
into $\kappa$ nonoverlapping subintervals,
and $p_{i}$ is the probability that failure time lies in
the $i$th subinterval.

In fact, we chose $\tilde{T}_{s}=
788400$, $946080$ and $998640$ and for each of them divided
the interval into subintervals shown in Table 3.
The values of the test statistic and
other variables in (\ref{eq:T}) are, for example,
as in Table 4 for the data in Table 2 and (\ref{eq:trueVal}).
%
%
\begin{table}[t]
\renewcommand{\thetable}{\arabic{table}}
\caption{Nonoverlapping subintervals}
\label{tab:group}
\begin{center}
\begin{tabular}{c|c|c|c|c}
  & \multicolumn{4}{|c}{subinterval}
  \\
  \cline{2-5}
  \raisebox{1.5ex}[0pt]{$\tilde{T}_{s}$}
  & 1 & 2 & 3 & 4
  \\
  \hline
  788400
  & $[0,12]$ & $(12,14]$ & $(14,\infty)$
  \\
  946080
  & $[0,16]$ & $(16,18]$ & $(18,\infty)$
  \\
  998640
  & $[0,11]$ & $(11,12]$ & $(12,13]$ & $(13,\infty)$
  \\  
\end{tabular}
\end{center}
\end{table}
%

%
%
\begin{table}[t]
\renewcommand{\thetable}{\arabic{table}}
\caption{Values of the test statistic}
\label{tab:test}
\begin{center}
\begin{tabular}{c|cccc|cccc|c}
  $\tilde{T}_{s}$
  & $m_{1}$ & $m_{2}$ & $m_{3}$ & $m_{4}$
  & $p_{1}$ & $p_{2}$ & $p_{3}$ & $p_{4}$
  & $T$
  \\
  \hline
  788400
  & $3$ & $4$ & $2$ &
  & $0.122186$ & $0.067105$ & $0.810709$ &
  & $26.22508$ 
  \\
  946080
  & $3$ & $3$ & $3$ &
  & $0.479264$ & $0.125719$ & $0.395017$ &
  & $3.572318$
  \\
  998640
  & $3$ & $3$ & $2$ & $3$
  & $0.225639$ & $0.059853$ & $0.064528$ & $0.649980$
  & $13.19004$ 
  \\  
\end{tabular}
\end{center}
\end{table}
%


In the simulation process
we sought parameter estimates,
calculated $T$ on simulated data,
and checked whether the value, say $T_{\rm sim}$, was at least as large as
the value of $T$ in Table 4.

\subsubsection*{Simulation result}
The result is shown in Table 5.
The second row indicates the number
of the data sets on that $T_{\rm sim}\geq T$ held.
The last column indicates the number of
the data sets 
on that the inequality held
simultaneously in the three $\tilde{T}_{s}$'s.
From this, we can reject the hypothesis
at any level $\alpha>1 \times 10^{-3}$ because
the p-value is less than $1 \times 10^{-3}$.

%
\begin{table}[t]
\renewcommand{\thetable}{\arabic{table}}
\caption{Successful number in satisfying $T_{sim}\geq T$}
\label{tab:res}
\begin{center}
\begin{tabular}{c|cccc}
  $\tilde{T}_{s}$ & $788400$ & $946080$ & $998640$ & simul.
  \\\hline
  $T_{\rm sim}\geq T$ & $1$  & $150$    & $16$     & $0$
\end{tabular}
\end{center}
\end{table}
%


The unsuccessful number was $158$
in finding a vector of ML estimates.
Both bias and variance of $\zeta$ were relatively large
as shown in Table 6.

%
\begin{table}[t]
\renewcommand{\thetable}{\arabic{table}}
\caption{Biases and variances}
\label{tab:bias}
\begin{center}
\begin{tabular}{c|cccc}
           & $\beta$    & $n$          & $\zeta$    & $\tilde{V}_{th}$
  \\\hline
  Bias     & $0.201057$ & $0.089020$ & $0.833870$ & $-0.053946$
  \\
  Variance & $0.496762$ & $0.114639$ & $7.444595$ & $0.016145$
\end{tabular}
\end{center}
\end{table}

%
%

%

%
\section{Conclusions}
\label{sec:conc}

Under the condition
that specimens were subjected to a normal level of stress and
did not fail before the test,
we considered the two-parameter Weibull CE model
with the threshold parameter
in the multiple-step SSALT. 
This consideration revealed that
the shape parameter $\beta$ must be larger than $1$ for
the model to fit the realistic situation in which the durability
of specimens decreases as they are used longer or
with higher stress.

After simplifying the model without loss of generality,
for a various sets of parameter values
we showed the average and standard deviation of failure
time versus the duration $\tilde{T}_{s}$ of the specimen's
use before the test.
As we have seen in Section \ref{sec:example},
we can utilize them in calculations for ML estimates.

In the section
we performed a goodness of fit test
by means of Monte Carlo simulation.
In general it is not clear whether the inverse power law
holds for every stress appearing in the step-stress test,
and it is not clear even whether the basic idea of the CE model
is available.
Including  these things,
the example we showed gives a way of checking
the validity of the CE model.

%

\vskip 3mm
\noindent {\bf Acknowledgments}
\vskip 3mm

The author would like to thank the referees for their
helpful comments to improve this paper.
%
\appendix
%
{\flushleft{\Large\bf Appendix}}:
{\large\bf Derivatives needed for the ML estimation}

\vspace{5mm}
{\noindent For}
$\theta_{1},\ \theta_{2}\in\{\beta,n,\zeta,\tilde{V}_{th}\}$, the
derivatives of the expressions in the left-hand side of
(\ref{eq:likeliEQ}) are given as follows:
\begin{eqnarray*}
  \lefteqn{\frac{\partial}{\partial\theta_{2}}
  \left\{
    \sum_{j=1}^{N}\frac{\lambda_{\theta_{1}}^{(j)}}{d_{j}}
    +\sum_{j=1}^{N}\delta_{\theta_{1}}(t_{1}^{(j)};\tilde{T}_{s}^{(j)})
  \right\}}\\
  & & = \sum_{j=1}^{N}
  \left\{
    -\left(\frac{1}{d_{j}}\frac{\partial
        d_{j}}{\partial\theta_{2}}\right)
    \frac{\lambda_{\theta_{1}}^{(j)}}{d_{j}}
    +\frac{1}{d_{j}}\frac{\partial\lambda_{\theta_{1}}^{(j)}}
    {\partial\theta_{2}}
    +\frac{\partial}{\partial\theta_{2}}\delta_{\theta_{1}}(t_{1}^{(j)};\tilde{T}_{s}^{(j)})
  \right\}\\
  & & = \sum_{j=1}^{N}
  \left\{
    -C_{\theta_{2}}\frac{\lambda_{\theta_{1}}^{(j)}}{d_{j}}
    \frac{\lambda_{\theta_{2}}^{(j)}}{d_{j}}
    +\frac{1}{d_{j}}\frac{\partial\lambda_{\theta_{1}}^{(j)}}
    {\partial\theta_{2}}
    +\frac{\partial}{\partial\theta_{2}}\delta_{\theta_{1}}(t_{1}^{(j)};\tilde{T}_{s}^{(j)})
  \right\},
\end{eqnarray*}
where
\begin{eqnarray*}
  C_{\theta}
  &\defeq&
  \left\{
    \begin{array}{rl}
      1 & (\theta=\beta),\\
      \frac{\beta}{\zeta} & (\theta=n),\\
      -\frac{\beta}{\zeta} & (\theta=\zeta),\\
      -\frac{\beta n}{\zeta} & (\theta=\tilde{V}_{th}),
    \end{array}
  \right.
  \quad E_{j}\defeq
  \exp\left(-\varepsilon^{\beta}(t_{l}^{(j)};\tilde{T}_{s}^{(j)})
      +\varepsilon^{\beta}(t_{l-1}^{(j)};\tilde{T}_{s}^{(j)})\right),
  \\
  \frac{\lambda_{\theta}^{(j)}}{d_{j}}
  &=&
  \frac{
    -\delta_{\theta}(t_{l-1}^{(j)};\tilde{T}_{s}^{(j)})
    +\delta_{\theta}(t_{l}^{(j)};\tilde{T}_{s}^{(j)})E_{j}
    }{1-E_{j}},\\
  \frac{1}{d_{j}}\frac{\partial\lambda_{\theta_{1}}^{(j)}}{\partial\theta_{2}}
  &=&
  \left[
    -\left\{
      \frac{\partial}{\partial\theta_{2}}
      \delta_{\theta_{1}}(t_{l-1}^{(j)};\tilde{T}_{s}^{(j)})
      -C_{\theta_{2}}\delta_{\theta_{1}}(t_{l-1}^{(j)};\tilde{T}_{s}^{(j)})
      \delta_{\theta_{2}}(t_{l-1}^{(j)};\tilde{T}_{s}^{(j)})
    \right\}
  \right.\\
  &&
  \left.
    \makebox[0.5mm]{}+\left\{
      \frac{\partial}{\partial\theta_{2}}
      \delta_{\theta_{1}}(t_{l}^{(j)};\tilde{T}_{s}^{(j)})
      -C_{\theta_{2}}\delta_{\theta_{1}}(t_{l}^{(j)};\tilde{T}_{s}^{(j)})
      \delta_{\theta_{2}}(t_{l}^{(j)};\tilde{T}_{s}^{(j)})
    \right\}E_{j}
  \right]\left/(1-E_{j}).\raisebox{0pt}[20pt][20pt]{}\right.
  \end{eqnarray*}
  In the above, note that the following equation does not necessarily
  hold because (\ref{eq:likeliEQ}) is a simplified equation not
  the original likelihood equation:
  \[
  \frac{\partial}{\partial\theta_{2}}
  \left\{\sum_{j=1}^{N}\frac{\lambda_{\theta_{1}}^{(j)}}{d_{j}}
    +\sum_{j=1}^{N}\delta_{\theta_{1}}(t_{1}^{(j)};\tilde{T}_{s}^{(j)})\right\}
  = \frac{\partial}{\partial\theta_{1}}
  \left\{\sum_{j=1}^{N}\frac{\lambda_{\theta_{2}}^{(j)}}{d_{j}}
    +\sum_{j=1}^{N}\delta_{\theta_{2}}(t_{1}^{(j)};\tilde{T}_{s}^{(j)})
  \right\}.
  \]
  Each derivative of $\delta_{\theta}\ (\theta\in \{\beta,n\zeta,
  \tilde{V}_{th}\})$ is given in the following. The arguments
  in some derivatives are omitted as far as it does not cause a confusion.
  \begin{eqnarray*}
    &&
    \dstyle{\frac{\partial\delta_{\beta}}{\partial\beta}=\delta_{\beta}
      \ln\varepsilon},\quad
    \dstyle{\frac{\partial\delta_{\beta}}{\partial n}=
      \left(C_{n}\ln\varepsilon+\frac{1}{\zeta}\right)\delta_{n}},\quad
    \dstyle{\frac{\partial\delta_{\beta}}{\partial \zeta}=
      \left(C_{\zeta}\ln\varepsilon-\frac{1}{\zeta}\right)\delta_{\zeta}},\\
    &&
    \dstyle{\frac{\partial\delta_{\beta}}{\partial \tilde{V}_{th}}=
      \left(C_{\tilde{V}_{th}}\ln\varepsilon-\frac{n}{\zeta}\right)
      \delta_{\tilde{V}_{th}}},\quad
    \dstyle{\frac{\partial\delta_{n}}{\partial \beta}}=
    \delta_{n}\ln\varepsilon,\\
    &&
    \lefteqn{\frac{\partial}{\partial n}\delta_{n}(t_{i};\tilde{T}_{s})}\\
    &&
    \makebox[3em]{}
    =\left(C_{n}-\frac{1}{\zeta}\right)
    \varepsilon^{\beta-2}(t_{i};\tilde{T}_{s})
    \left\{
      \frac{1}{\tilde{K}_{0}}
      \sum_{m=k}^{i}((m-1)\del{\tilde{V}}-\tilde{V}_{th})^{n}
      \ln ((m-1)\del{\tilde{V}}-\tilde{V}_{th})
    \right\}^{2}\\
    &&
    \makebox[6.5em]{}
    +\varepsilon^{\beta-1}(t_{i};\tilde{T}_{s})
    \left\{
      \frac{1}{\tilde{K}_{0}}
      \sum_{m=k}^{i}((m-1)\del{\tilde{V}}-\tilde{V}_{th})^{n}
      \Bigl(\ln ((m-1)\del{\tilde{V}}-\tilde{V}_{th})\Bigr)^{2}
    \right\},\\
    &&
    \dstyle{\frac{\partial\delta_{n}}{\partial\zeta}}=
      \left(C_{\zeta}+\frac{1}{\zeta}\right)\delta_{n},\\
    &&
    \frac{\partial}{\partial \tilde{V}_{th}}
    \delta_{n}(t_{i};\tilde{T}_{s})
    =\left(C_{\tilde{V}_{th}}+\frac{n}{\zeta}\right)
    \varepsilon^{\beta-2}(t_{i};\tilde{T}_{s})
    \left\{
      \frac{1}{\tilde{K}_{0}}
      \sum_{m=k}^{i}((m-1)\del{\tilde{V}}-\tilde{V}_{th})^{n-1}
    \right\}\\
    &&
    \makebox[7.5em]{}
    \times\left\{
      \frac{1}{\tilde{K}_{0}}
      \sum_{m=k}^{i}((m-1)\del{\tilde{V}}-\tilde{V}_{th})^{n}
      \ln ((m-1)\del{\tilde{V}}-\tilde{V}_{th})
    \right\}-\delta_{\tilde{V}_{th}}(t_{i};\tilde{T}_{s})
    \\
    &&
    \makebox[6.5em]{}
    -n\varepsilon^{\beta-1}(t_{i};\tilde{T}_{s})
    \left\{
      \frac{1}{\tilde{K}_{0}}
      \sum_{m=k}^{i}((m-1)\del{\tilde{V}}-\tilde{V}_{th})^{n-1}
      \ln ((m-1)\del{\tilde{V}}-\tilde{V}_{th})
    \right\},
    \\
    &&
    \dstyle{\frac{\partial\delta_{\zeta}}{\partial \beta}}=
    \delta_{\beta},\quad
    \dstyle{\frac{\partial\delta_{\zeta}}{\partial n}}=
    C_{n}\delta_{n},\quad
    \dstyle{\frac{\partial\delta_{\zeta}}{\partial \zeta}}=
    C_{\zeta}\delta_{\zeta},\quad
    \dstyle{\frac{\partial\delta_{\zeta}}{\partial \tilde{V}_{th}}}=
    C_{\tilde{V}_{th}}\delta_{\tilde{V}_{th}},\quad
    \dstyle{\frac{\partial\delta_{\tilde{V}_{th}}}{\partial \beta}}=
    \delta_{\tilde{V}_{th}}\ln\varepsilon,
  \end{eqnarray*}
  \begin{eqnarray*}
    &&
    \lefteqn{\frac{\partial}{\partial n}
      \delta_{\tilde{V}_{th}}(t_{i};\tilde{T}_{s})}\\
    &&
    \makebox[3em]{}
    =\left(C_{n}-\frac{1}{\zeta}\right)
    \varepsilon^{\beta-2}(t_{i};\tilde{T}_{s})
    \left\{
      \frac{1}{\tilde{K}_{0}}
      \sum_{m=k}^{i}((m-1)\del{\tilde{V}}-\tilde{V}_{th})^{n}
      \ln((m-1)\del{\tilde{V}}-\tilde{V}_{th})
    \right\}\\
    &&
    \makebox[6.5em]{}
    \times\left\{
      \frac{1}{\tilde{K}_{0}}
      \sum_{m=k}^{i}((m-1)\del{\tilde{V}}-\tilde{V}_{th})^{n-1}
    \right\}\\
    &&
    \makebox[4em]{}
    +\varepsilon^{\beta-1}(t_{i};\tilde{T}_{s})
    \left\{
      \frac{1}{\tilde{K}_{0}}
      \sum_{m=k}^{i}((m-1)\del{\tilde{V}}-\tilde{V}_{th})^{n-1}
      \ln ((m-1)\del{\tilde{V}}-\tilde{V}_{th})
    \right\},
    \\
    &&
    \dstyle{\frac{\partial\delta_{\tilde{V}_{th}}}{\partial \zeta}}=
    \left(C_{\zeta}+\frac{1}{\zeta}\right)\delta_{\tilde{V}_{th}},
    \\
    &&
    \frac{\partial}{\partial \tilde{V}_{th}}
    \delta_{\tilde{V}_{th}}(t_{i};\tilde{T}_{s})
    =\left(C_{\tilde{V}_{th}}+\frac{n}{\zeta}\right)
    \varepsilon^{\beta-2}(t_{i};\tilde{T}_{s})
    \left\{
      \frac{1}{\tilde{K}_{0}}
      \sum_{m=k}^{i}((m-1)\del{\tilde{V}}-\tilde{V}_{th})^{n-1}
    \right\}^{2}\\
    &&
    \makebox[7em]{}
    -(n-1)\varepsilon^{\beta-1}(t_{i};\tilde{T}_{s})
    \left\{
      \frac{1}{\tilde{K}_{0}}
      \sum_{m=k}^{i}((m-1)\del{\tilde{V}}-\tilde{V}_{th})^{n-2}
    \right\}.
  \end{eqnarray*}

%
%
\vskip 3mm
\clearpage
\noindent {\bf REFERENCES}
\vskip 3mm

\noindent Bagdanavicius, V.B. (1978)
A statistical test of a model of additive accumulation of damage,
{\it Theory of Probability and its Applications},
23 (2), pp. 385--390.
\vskip 3mm

\noindent Bank, R.E. \& Rose, D.J. (1981)
Global approximate Newton methods,
{\it Numerische Mathematik},
37, pp. 279--295.
\vskip 3mm

\noindent Hirose, H. (1996)
Theoretical foundation for residual lifetime estimation,
{\it Transactions of the Institute of Electrical Engineers of Japan},
116-B (2), pp. 168--173.
\vskip 3mm

\noindent Hirose, H. (1997)
Mixture model of the power law,
{\it IEEE Transactions on Reliability},
46 (1), pp. 146--153.
\vskip 3mm

\noindent Khamis, I.H. (1997)
Optimum $m$-step, step-stress design with $k$ stress variables,
{\it Communications in Statistics---Simulation and Computation},
26 (4), pp. 1301--1313.
\vskip 3mm

\noindent Lu, M.-W. \& Rudy, R.J. (2002)
Step-stress accelerated test,
{\it International Journal of Materials \& Product Technology},
17 (5-6), pp. 425--434.
\vskip 3mm

\noindent McSorley, E.O., Lu, J.-C. \& Li, C.-S. (2002)
Performance of parameter-estimates in step-stress accelerated
life-tests with various sample-sizes,
{\it IEEE Transactions on Reliability},
51 (3), pp. 271--277.
\vskip 3mm

\noindent Miller, R. \& Nelson, W. (1983)
Optimum simple step-stress plans for accelerated life testing,
{\it IEEE Transactions on Reliability},
R-32 (1), pp. 59--65.
\vskip 3mm

\noindent Miner, M.A. (1945)
Cumulative damage in fatigue,
{\it Journal of Applied Mechanics},
12, pp. 159--164.
\vskip 3mm

\noindent Nelson, W. (1980)
Accelerated life testing -- step-stress models and data analyses,
{\it IEEE Transactions on Reliability},
R-29 (2), pp. 103--108.
\vskip 3mm

\noindent Nelson, W. (1990)
{\it Accelerated Testing: Statistical Models, Test Plans, and Data Analyses}
(New York, Wiley).
\vskip 3mm

\noindent Park, S.-J.  \& Yum, B.-J. (1998)
Optimal design of accelerated life tests under modified stress loading
methods,
{\it Journal of Applied Statistics},
25 (1), pp. 41--62.
\vskip 3mm

\noindent Rao, C.R. (2002)
{\it Linear Statistical Inference and Its Applications}
(New York, Wiley).
\vskip 3mm

\noindent Ross, S.M. (2002)
{\it Simulation} (New York, Academic Press).
\vskip 3mm

\noindent Xiong, C. (1999)
Step stress model with threshold parameter, 
{\it Journal of Statistical Computation and Simulation},
63, pp. 349--360.
\vskip 3mm

\noindent Yeo, K.-P. \& Tang, L.-C. (1999)
Planning step-stress life-test with a target acceleration-factor,
{\it IEEE Transactions on Reliability},
48 (1), pp. 61--67.
%
%
%
%
\end{document}